\documentclass[draft]{amsart}
\usepackage{amssymb,amsxtra}
\usepackage[french,english]{babel}
\usepackage[T1]{fontenc}

\numberwithin{equation}{section}

\newtheorem{Theoreme}{Th\'eor\`eme}[section] 

\newtheorem{prop}[Theoreme]{Proposition} 
\newtheorem{Assertion}[Theoreme]{Assertion} 
\newtheorem{Lemme}[Theoreme]{Lemme} 
\newtheorem{cor}[Theoreme]{Corollaire} 

\theoremstyle{definition}
\newtheorem{defn}[Theoreme]{D\'efinition} 

\theoremstyle{remark}
\newtheorem{Remarque}[Theoreme]{Remarque} 
\newtheorem{Remarques}[Theoreme]{Remarques}

\newenvironment{dem}{\proof[D\'emonstration]}{\endproof}

\newcommand{\N}{{\mathbb{N}}}
\newcommand{\C}{{\mathbb{C}}}
\newcommand{\R}{{\mathbb{R}}}
\newcommand{\Z}{{\mathbb{Z}}}

\newcommand{\aalg}{\mathfrak{a}}
\newcommand{\fa}{\mathfrak{a}}

\newcommand{\fsl}{\mathfrak{sl}}

\newcommand{\calH}{\mathcal{H}}
\newcommand{\hec}{\mathcal{H}}
\newcommand{\calP}{\mathcal{P}}
\newcommand{\calS}{\mathcal{S}}
\newcommand{\ima}{{I}_\chi}

\newcommand{\V}{V_\chi}
\renewcommand{\v}{v_\chi}
\newcommand{\half}{{\frac{1}{2}}}

\newcommand{\CW}{{\mathbb{C} W}}

\def\Hom{\operatorname {Hom}}
\def\End{\operatorname {End}}
\def\boplus{\mathbin{\boldsymbol{\oplus}}}

\newcommand{\dual}[2]{{\langle {#1} ,  {#2} \rangle}}

\newcommand{\pscal}[2]{{( {#1}  \mid {#2} )}}

\newcommand{\sgn}{\mathrm{sgn}}
\newcommand{\triv}{\mathrm{triv}}
\newcommand{\pnr}[2]{P_{{#1},{#2}}}
\newcommand{\Pnr}{{P_{n,r}}}
\newcommand{\Pnpr}{{P_{n+1,r}}}  

\begin{document}

\selectlanguage{french}

\title[\phantom{.}]{Repr\'esentations de dimension finie de
  l'alg\`ebre de Cherednik rationnelle}
\author[C.~Dez\'el\'ee]{Charlotte Dez\'el\'ee}
\address{D\'epartement de math\'ematiques, Universit\'e de Brest,
  29285 Brest cedex, France} \date{16 janvier 2002}
\email{Charlotte.Dezelee@univ-brest.fr}
\begin{abstract} On donne une condition n\'ecessaire et suffisante pour
  l'existence de modules de dimension finie sur l'alg\`ebre de
  Cherednik rationnelle associ\'ee à un système de racines.
\end{abstract}

\maketitle


\section{Introduction et notations}
\label{sec1}

Soit $\aalg_\mathbb{R}$ un $\mathbb{R}$-espace vectoriel de
dimension $\ell \ge 1$, $R \subset \aalg_\mathbb{R}^*$ un
syst\`eme de racines réduit, $W$ le groupe de Weyl correspondant
et $\aalg = \mathbb{C} \otimes_\mathbb{R} \aalg_\mathbb{R}$. On
notera $r_{\alpha}$ la r\'eflexion associ\'ee \`a la racine
$\alpha$ et $\alpha\spcheck \in \fa$ la coracine de $\alpha$.  On
fixe une fonction de multiplicit\'e $k:R \to \C$ sur l'ensemble
des racines, donc $k_{w(\alpha)} = k_\alpha$ pour tous $\alpha \in
R$, $w \in W$.

Soit $\calP = S(\fa^*) = \boplus_{n \geq 0}\calP_n$, o\`u $\calP_n
=S^{n}(\aalg^{*})$, l'alg\`ebre sym\'etrique de $\fa^*$.  On pose
$\calP_+ = \boplus_{n \ge 1} \calP_n$.  Un \'el\'ement $f \in
\calP$ peut \^etre identifi\'e \`a la multiplication par $f$ dans
$\End_\C(\calP)$ et l'on fait op\'erer $W$ de fa\c{c}on naturelle
sur $\calP$.  Si $\partial_y$ est le champ de vecteurs associ\'e
\`a $y \in \fa$ on d\'efinit alors l'op\'erateur de Dunkl $T_y
=T_y(k) \in \End_\C(\calP)$ par
$$
T_y(k) = \partial_y + \half\sum_{\alpha \in R} k_\alpha
\frac{\dual{\alpha}{y}}{\alpha} (1 - r_\alpha)= \partial_y
+\sum_{\alpha \in R^+} k_\alpha \frac{\dual{\alpha}{y}}{\alpha} (1
- r_\alpha)
$$
où $R^+ \subset R$ est un système de racines positives.  On
sait que $T : y \mapsto T_y$ s'\'etend en un isomorphisme
d'alg\`ebres de $S(\fa)$ sur $\calS= \calS(k)=\C[T_y : y \in \fa]$
(cf., par exemple, \cite[Theorem~2.12]{DO}).  On posera $\calS_n =
T(S^n(\fa))$, de sorte que $\calS = \C \boplus \calS_+$ avec
$\calS_+ = \boplus_{n \ge 1} \calS_n$.

L'alg\`ebre de Cherednik rationnelle (cf.~\cite{EG}), not\'ee
$\hec(k)$ ou $\hec$, est la sous-alg\`ebre de $\End_\C(\calP)$
engendr\'ee par les $w \in W$, $x \in \aalg^*$ et les $T_y$, $y
\in \aalg$. Ces g\'en\'erateurs sont li\'es par les relations
suivantes:
\begin{enumerate}
\item $[T_{y},x] = \dual{y}{x} +\frac{1}{2}\sum_{\alpha \in R}
  k_{\alpha}\dual{y}{\alpha} \dual{\alpha\spcheck}{x}r_{\alpha}$;
\item $wxw^{-1}=w(x)$;
\item $wT_{y}w^{-1}=T_{w(y)}$.
\end{enumerate}
Rappelons \cite[Corollary~4.4]{EG} que $\hec$ v\'erifie un
th\'eor\`eme de Poincar\'e-Birkhoff-Witt (PBW). En effet, si l'on
pose pour tout $n \in \Z$
$$
\hec_n =\textstyle{\bigoplus_{j -i=n, w \in W}}
\mathbb{C}\calS_{i}w\calP_{j} = \textstyle{\bigoplus_{j-i=n, w \in
    W}} \mathbb{C}\calP_{j}w \calS_{i},
$$
on a alors
$$
\hec=\textstyle{\bigoplus_{n \in \Z}} \hec_n =\calP
\otimes\mathbb{C}W \otimes \calS = \calS \otimes\mathbb{C}W
\otimes \calP.
$$

Il r\'esulte de \cite[Theorem~2.2]{BEG} que pour des valeurs
g\'en\'eriques de la fonction $k$, l'alg\`ebre $\calH(k)$ ne
poss\`ede pas de repr\'esentation de dimension finie.  L'objet de
ce travail est de chercher des conditions n\'ecessaires et
suffisantes pour qu'il existe des $\hec$-modules de dimension
finie et de d\'eterminer certaines de leurs propri\'et\'es.  Les
principaux r\'esultats obtenus sont les suivants. Dans la
section~\ref{sec4} on montre (Th\'eor\`eme~\ref{thm41}) que tout
$\hec(k)$-module irr\'eductible de dimension finie est isomorphe
\`a l'unique quotient simple d'un module de Verma g\'en\'eralis\'e
(cf.~\cite{DO}); on en d\'eduit \`a la section~\ref{sec5} une
caract\'erisation des multiplicit\'es $k$ pour lesquelles des
$\hec(k)$-modules de dimension finie existent et une description
de ces derniers (Th\'eor\`eme~\ref{thm54}, Remarque~\ref{rem55} et
Proposition~\ref{prop56}).  La section~\ref{sec6} est consacr\'ee
\`a des exemples, notamment le cas d'un groupe de Weyl en rang $2$
y est (presque) compl\`etement trait\'e.


\section{Propri\'et\'es de $\calH$}
\label{sec2}

On d\'efinit \cite[p.~144]{Bou} une forme bilin\'eaire
sym\'etrique d\'efinie positive $W$-in\-va\-riante sur
$\fa_\R^{*}$ (que l'on étend à $\fa^*$) en posant
$B^{*}(x,z)=\sum_{\alpha \in R}
\dual{\alpha\spcheck}{x}\dual{\alpha\spcheck}{z}$.  On peut alors
identifier $\aalg^{*}$ \`a $\aalg$ via $x \mapsto B(x)$, o\`u
$B(x)$ est caract\'eris\'e par $B^{*}(x,z)=\dual{B(x)}{z}$.  Ainsi
$B$ est un isomorphisme $W$-lin\'eaire et l'on a: $B(\alpha)=
\frac{B^{*}(\alpha,\alpha)}{2} \alpha\spcheck$,
$\dual{B(x)}{B^{-1}(y)}=\dual{y}{x}$, pour tous $\alpha \in R$ et
$(x,y) \in \aalg^{*} \times \aalg$.  On en d\'eduit que:
$$
\sum_{\alpha \in R}
k_{\alpha}\dual{y}{\alpha}\dual{\alpha\spcheck}{x}r_{\alpha} =
\sum_{\alpha \in R}
k_{\alpha}\dual{B(x)}{\alpha}\dual{\alpha\spcheck}{B^{-1}(y)}r_{\alpha}.
$$
En utilisant cette relation on montre que l'on peut d\'efinir
un anti-automorphisme involutif $\sigma$ de $\hec$ par:
$$
\sigma(x)=T_{B(x)}, \quad \sigma(T_y)= B^{-1}(y), \quad
\sigma(w)= w^{-1},
$$
pour tous $x \in \aalg^{*}, y \in \aalg, w \in W$.
On peut aussi d\'efinir un automorphisme $\phi$ de $\hec$ (d'ordre
$4$) en posant:
$$
\phi(x)=-T_{B(x)}, \quad \phi(T_y)= B^{-1}(y), \quad \phi(w)=
w.
$$
On remarquera que $\phi$ et $\sigma$ \'echangent $\hec_n$ et
$\hec_{-n}$ pour tout $n\in\mathbb{Z}$.

\subsection*{Une application bilin\'eaire sur $\hec$}
Gr\^ace \`a PBW il vient: $\hec=\mathbb{C}W \boplus (\calP_{+}
\hec + \hec \calS_{+})$.  On peut donc d\'efinir la projection
$\pi : \hec \rightarrow \mathbb{C}W $ parall\`element \`a
$\calP_{+} \hec + \hec \calS_{+}$.  Comme $W$ laisse stables
$\calP_{+}$ et $\calS_{+}$, $\pi$ est un morphisme de $W$-modules
pour l'action par multiplication \`a gauche, ou \`a droite de $W$,
i.e.~$\pi(w h) = w\pi(h)$ et $\pi(hw) = \pi(h)w$, $h \in \calH$,
$w \in W$; observons \'egalement que $\pi(\sigma(h)) =
\sigma(\pi(h))$.  D\'efinissons une application bilin\'eaire
$\beta : \hec \times \hec \to \C W$ par:
$$
\forall a,b \in \hec, \quad \beta(a,b)= \pi(\sigma(a)b).
$$

Les assertions du lemme suivant r\'esultent de calculs
imm\'ediats.

\begin{Lemme} \label{lem21}
  On a, pour tous $a,b,h \in \hec $ et $w,s \in W$:

\noindent {\rm 1.} $\beta(a,b):=\sigma(\beta(b,a))$;

\noindent {\rm 2.}  $R(\beta)=\{a\in\hec : \forall b \in\hec, \,
\beta(a,b)=0\}=\{a\in\hec : \forall b \in\hec, \, \beta(b,a)=0\}$;

\noindent {\rm 3.}  $\beta(aw,bs)=w^{-1}\beta(a,b)s$ (on dira que 
$\beta$ est $\sigma$-lin\'eaire);

\noindent {\rm 4.} $\beta(ah,b)=\beta(a,\sigma(h)b)$ (on dira que
$\beta$ est $\sigma$-sym\'etrique);

\noindent {\rm 5.} $R(\beta)$ contient $\hec \calS_{+}$ et
$\sigma(\hec \calS_{+})=\calP_{+} \hec $.
\end{Lemme}

Il d\'ecoule de PBW que $\calH_p \subset \calP_+ \calH + \calH
\calS_+$ lorsque $p \ne0$; on en d\'eduit
$$
\beta(\calH_{m},\calH_n) = \pi(\sigma(\calH_{m})\calH_n) =
\pi(\calH_{-m} \calH_n) \subset \pi(\calH_{n-m}) = 0
$$
pour tous $m \ne n$.

Comme $\hec= (\calP\otimes \mathbb{C}W)\boplus\hec \calS_{+}$, on
peut consid\'erer $\beta$ comme une application bilin\'eaire
$\sigma$-sym\'etrique sur $\hec / \hec \calS_{+} \simeq
\calP\otimes \mathbb{C}W$ en posant:
$$
\forall p,q \in \calP, \ \, \forall w,s \in \mathbb{C}W, \quad
\beta(p \otimes w,q \otimes s)= \sigma(w)\beta(p,q)s.
$$
On remarquera que $\beta(\calP_n\otimes
\mathbb{C}W,\calP_m\otimes \mathbb{C}W)=0$ si $n \neq m$.  Donc
$\beta$ est d\'etermin\'ee par ses restrictions aux espaces
vectoriels de dimension finie $\calP_n \otimes \mathbb{C}W$.

\subsection*{\'Equivalence de cat\'egories}
Soit $\tau : W \to \{\pm 1\}$ une repr\'esentation de
dimension~$1$.  D\'efinissons une fonction de multiplicit\'e
$k^\tau : R \to \C$ par $k^\tau(\alpha) = \tau(r_\alpha)
k_\alpha$.  Il r\'esulte des relations 1,2,3 entre les
g\'en\'erateurs de $\calH(k)$ rappel\'ees au \S~{\ref{sec1}} que
l'on peut d\'efinir un morphisme d'alg\`ebre $\epsilon_\tau :
\calH(k) \to \calH(k^\tau)$ en posant
$$
\epsilon_\tau(x) = x, \quad \epsilon_\tau(T_y(k)) =
T_y(k^\tau), \quad \epsilon_\tau(w) = \tau(w) w
$$
pour tous $x \in \fa^*$, $y \in \fa$, $w \in W$. Il est clair
que $\epsilon_\tau$ est un isomorphisme dont on notera encore
$\epsilon_\tau$ l'inverse. Si $M$ est un $\calH(k^\tau)$-module,
on peut alors d\'efinir un $\calH(k)$-module $M^{\epsilon_\tau}$
en munissant le $\C$-espace vectoriel de l'action $h.v =
\epsilon_\tau(h)v$ pour tous $h \in \calH(k)$, $v \in M$. On en
d\'eduit ainsi une \'equivalence de cat\'egories, $M \mapsto
M^{\epsilon_\tau}$, entre $\calH(k^\tau)$-$\mathrm{mod}$ et
$\calH(k)$-$\mathrm{mod}$ (on a
$\Hom(M^{\epsilon_\tau},N^{\epsilon_\tau}) = \Hom(M,N)$).

\begin{Remarque} \label{rem22}
  Si $V$ est un $W$-module on notera \'egalement
  $V^{\epsilon_\tau}$ le $W$-module obtenu en munissant le
  $\C$-espace vectoriel $V$ de l'action $w.v=\epsilon_\tau(w)v$.
  Le $W$-module $V^{\epsilon_\tau}$ est alors isomorphe \`a
  $V\otimes_{\mathbb{C}W}V_{\tau}$, o\`u $V_{\tau}$ d\'esigne un
  $W$-module irr\'eductible de type $\tau$; en particulier si
  $V_{\chi}$ est un $W$-module irr\'eductible de type $\chi$,
  alors $V_{\chi}^{\epsilon_\tau}$ est un $W$-module
  irr\'eductible de type $\chi\otimes\tau$.
\end{Remarque}

On note $\sgn$ le caract\`ere $w \mapsto \det(w)$ de $W \subset
\mathrm{GL}(\fa)$. La construction pr\'ec\'edente s'applique \`a
$\tau=\sgn$ et (pour simplifier) on posera $\varepsilon =
\epsilon_\sgn$. Le foncteur $M \mapsto M^\varepsilon$ \'etablit
donc une \'equivalence de cat\'egories entre
$\calH(-k)$-$\mathrm{mod}$ et $\calH(k)$-$\mathrm{mod}$.


\section{Modules de Verma sur $\calH$}
\label{sec3}

On note $W\sphat \,$ l'ensemble des caract\`eres irr\'eductibles
de $W$. Tout $\calH$-module $M$ \'etant un $W$-module, il se
d\'ecompose en $M =\boplus M[\chi]$ o\`u $M[\chi]$ est la
composante isotypique de type $\chi$ de $M$.  Soit $\chi \in
W\sphat\,$ et $V_\chi$ un $W$-module irr\'eductible de type
$\chi$.  Rappelons la d\'efinition d'un module de Verma de plus
bas poids $\chi$ introduite dans \cite[(25)]{DO}.  On munit
$V_\chi$ d'une structure de $\calS \otimes \C W$-module en posant
$\calS_+.V_\chi =0$.

\begin{defn}
\label{def31}
On appelle module de Verma de plus bas poids $\chi$, not\'e
$M(\chi)=M(\chi,k)$, le module induit par $V_\chi$ de $\calS
\otimes \CW$ \`a $\hec$:
$$
M(\chi)= \mathrm{ind}_{\calS\otimes
  \mathbb{C}W}^{\hec}(V_\chi)= \hec \otimes_{\calS \otimes
  \mathbb{C}W} V_\chi.
$$
\end{defn}

Il r\'esulte de $\calH = \calP \otimes \CW \boplus \calH\calS_+$
que $M(\chi)$ s'identifie \`a $\calP\otimes \V$ comme
$\calP$-module.  Les propri\'et\'es \'enonc\'ees dans la
proposition qui suit d\'ecoulent de \cite[2.5]{DO}. Rappelons que
$\calP$ poss\`ede une structure naturelle de $\calH$-module.

\begin{prop} 
\label{prop32}
{\em (a)} Si $\chi = \triv$ est le caract\`ere trivial, $M(\triv)$
s'identifie au $\calH$-module $\calP$.

\noindent {\em (b)} Si $\v \in \V \smallsetminus \{0\}$ et $\ima$ est
l'annulateur de $\v$ dans $\mathbb{C}W$, on a $M(\chi) \simeq
\hec.\v \simeq \hec / (\hec \calS^{+}+\hec\ima)$.

\noindent  {\em (c)} Tout sous-module $M$ de $M(\chi)$ est gradu\'e:
$M=\boplus_{n\geq 0}M_n$ o\`u $M_n = (\calP_n \otimes\V) \cap M$.

\noindent {\em (d)} Un sous-module $M$ de $M(\chi)$ est propre si
et seulement si $M\cap\V=0$.

\noindent {\em (e)} $M(\chi)$ admet un unique sous-module maximal,
et donc un unique quotient simple que l'on note $L(\chi) =
L(\chi,k)$.

\noindent {\em (f)} Soit $V$ un $\calH$-module engendr\'e par $v$
tel que $\CW.v \simeq V_\chi$ et $\calS_+.v= 0$. Il existe alors
un morphisme surjectif de $\calH$-modules $M(\chi) \to V$; si $V$
est irréductible on a $V \simeq L(\chi)$.  \qed
\end{prop}

Signalons le corollaire:

\begin{cor}
\label{cor33}
Soit $\tau$ une repr\'esentation de dimension~$1$ de $W$.  Il
existe un isomorphisme naturel de $\hec(k)$-modules
$$
M(\chi,k) \simeq M(\chi \otimes \tau, k^\tau)^{\epsilon_\tau}.
$$
\end{cor}

\begin{dem}
  On applique le~(b) la proposition précédente, dont on adopte les
  notations.  Remarquons que si $J$ est un id\'eal \`a gauche de
  $\hec(k^\tau)$ et $M=\hec(k^\tau)/J$, alors $M^{\epsilon_\tau}$
  est isomorphe au $\hec(k)$-module $\hec(k)/\epsilon_\tau(J)$. Le
  corollaire découle de cette remarque appliquée à
  $J=\hec(k^\tau)\calS_{+}(k^\tau) +
  \hec(k^\tau)I_{\chi\otimes\tau}$.  En effet, fixons l'annulateur
  $I_{\chi \otimes \tau}$ d'un \'el\'ement non nul de $V_{\chi
    \otimes \tau}$; alors, $I_\chi=\epsilon_\tau(I_{\chi \otimes
    \tau})$ est l'annulateur de ce m\^eme \'el\'ement dans $V_\chi
  = V_{\chi \otimes \tau}^{\epsilon_\tau}$.  Donc
  $\epsilon_\tau(J) = \hec(k)\calS_{+}(k)+\hec(k)I_{\chi}$, d'où
  le résultat voulu.
\end{dem}

\subsection*{Propri\'et\'es de $L(\chi)$}
Comme il est remarqu\'e en \cite[2.6]{DO} on peut munir $M(\chi)$
d'une forme analogue \`a la forme de Shapovalov. Nous donnons
ci-dessous une mani\`ere de construire une telle forme, ce qui
nous servira au \S~\ref{sec5}.

D\'efinissons tout d'abord une forme bilin\'eaire sur $\hec$ de la
fa\c{c}on suivante.  On fixe $0 \ne \v \in \V$ et une forme
bilin\'eaire sym\'etrique non d\'eg\'en\'er\'ee $W$-invariante
$\pscal{\phantom{.}}{\phantom{.}}$ sur $\V$ (on peut la prendre
sym\'etrique puisque $\V\simeq\V^{*}$ comme $W$-module).  On a
donc, pour tous $u,v\in \V$ et $w \in \mathbb{C}W$,
$\pscal{w.u}{v}=\pscal{u}{\sigma(w).v}$.  On peut alors d\'efinir
une forme bilin\'eaire sur $\hec$ en posant :
$$
\pscal{a}{b}_\chi =\pscal{\v}{\beta(a,b).\v}.
$$
Remarquons que $\pscal{\phantom{.}}{\phantom{.}}_\chi$ est non
nulle (sinon $\pscal{\v}{\V} = 0$).  On d\'eduit facilement des
propri\'et\'es de $\beta$ que
$\pscal{\phantom{.}}{\phantom{.}}_\chi$ est sym\'etrique et que
son radical contient $\hec S^{+} + \hec \ima$.

Comme $M(\chi)=\hec.\v \simeq \hec / (\hec S^{+} + \hec \ima)$, la
forme $\pscal{\phantom{.}}{\phantom{.}}_\chi$ induit une forme
bilin\'eaire sym\'etrique non nulle sur $M(\chi)$ par la formule:
$$
(a.\v,b.\v) = \pscal{a}{b}_\chi=\pscal{\v}{\beta(a,b).\v}
$$
pour tous $a,b \in \calH$.  Cette forme d\'epend des choix de
$\pscal{\phantom{.}}{\phantom{.}}$ et $\v$, nous allons voir que
son radical $R(\chi)$ n'en d\'epend pas.

\begin{Remarques} \label{rem34}
  (1) En utilisant la $\sigma$-sym\'etrie de $\beta$, on montre
  celle de $(\phantom{.},\phantom{.})$, qui est en particulier
  $W$-invariante ; il en r\'esulte que le radical $R(\chi)$ est un
  sous-$\hec$-module de $M(\chi)$, diff\'erent de $M(\chi)$.

\noindent (2)  Soit $M(\chi)=\boplus_{\tau \in W\sphat} \,
M(\chi)[\tau]$ la d\'ecomposition en composantes isotypiques du
$W$-module $M(\chi)$. On montre facilement que
$(M(\chi)[\tau],M(\chi)[\psi])=0$ si $\tau \neq \psi$.

\noindent (3) En identifiant $M(\chi)$ et $\calP\otimes \V$, on
d\'eduit facilement de $\beta(\calP_n\otimes \V,\calP_m\otimes
\V)=0$ pour $n \neq m$ que $(M_n(\chi),M_m(\chi))=0$ si $n \neq
m$.  Ou encore, du fait que $\beta(\hec_n,\hec_m)=0$ si $n\neq m$,
il d\'ecoule que $(\hec_n.\v,\hec_m.\v)=0$ si $n\neq m$.  En
observant que chaque $M_n(\chi)$ est un $W$-module, le (2)
implique alors $(M_n(\chi)[\tau],M_m(\chi)[\psi])=0$ si $n \ne m$
ou $\tau \ne \psi$.
\end{Remarques}

\begin{prop}
\label{prop35}
{\em 1.} Le radical $R(\chi)$ est l'unique sous-module maximal de
$M(\chi)$. Par cons\'equent $M(\chi)/R(\chi)$ est l'unique
quotient simple $L(\chi)$ de $M(\chi)$.

\noindent {\em 2.} Il existe un isomorphisme naturel $L(\chi,k)
\simeq L(\chi \otimes \tau,k^\tau)^{\epsilon_\tau}$.
\end{prop}

\begin{dem}
  1.  Soit $M$ un sous-$\hec$-module propre de $M(\chi)$. Alors
  $M$ \'etant gradu\'e et tel que $M\cap\V=0$, on a
  $M=\boplus_{n>0}M_n$. Par la Remarque~\ref{rem34}(3) et le fait
  que $\v \in M_{0}(\chi)$ il vient $(\v,M)=0$.  D'o\`u, par
  $\sigma$-sym\'etrie, $(\hec.\v,M) =(\v,M)=0$.  Donc $M$ est
  inclus dans $R(\chi)$.
  
  2. Par l'\'equivalence de cat\'egories entre
  $\hec(k)$-$\mathrm{mod}$ et $\hec(k^\tau)$-$\mathrm{mod}$, le
  $\hec(k)$-module $M(\chi\otimes \tau,k^\tau)^{\epsilon_\tau}$
  admet un unique quotient simple $L(\chi\otimes
  \tau,k^\tau)^{\epsilon_\tau}$. Mais $M(\chi,k)\simeq
  M(\chi\otimes \tau,k^\tau)^{\epsilon_\tau}$, par cons\'equent
  les $\hec(k)$-modules $L(\chi\otimes
  \tau,k^\tau)^{\epsilon_\tau}$ et $L(\chi,k)$ sont isomorphes.
\end{dem}

\begin{Remarque}
\label{rem36}
Lorsque $\chi$ est le caract\`ere trivial $R(\chi)$ co\"{\i}ncide
avec le radical $R(k)$ de la forme $(\phantom{.},\phantom{.})_k$
d\'efinie dans \cite{DOJ}. Ceci résulte des
Propositions~\ref{prop35} et~\ref{prop32}(e), et du fait que
$\calP/R(k) \simeq L(\triv,k)$ (cf.~\cite[2.6]{DO}). Si $\tau$ est
une repr\'esentation de dimension~$1$ de $W$ on a un isomorphisme
$L(\tau,k)\simeq L(\triv,k^\tau)^{\epsilon_\tau}$; en particulier,
$L(\sgn, k)\simeq L(\triv,-k)^{\varepsilon}$
\end{Remarque}


\section{$\hec$-modules irr\'eductibles de dimension finie}
\label{sec4}

Le but de ce paragraphe est de montrer que tout $\hec$-module
irr\'eductible de dimension finie est isomorphe \`a un $L(\chi)$
pour un $\chi \in W\sphat$.  Pour ce faire on va utiliser une
copie de $\fsl(2,\C)$ contenue dans l'algèbre $\hec^W$ des
$W$-invariants. Soit $\{z_1,\dots,z_\ell\}$ une base orthonormée
de $\fa_\R^*$ et $e_j = B(z_j)$, $1 \le j \le \ell$. Posons
$$
E=\frac{1}{2}\sum_{i=1}^{\ell}{z_i}^{2}, \quad
F=-\frac{1}{2}\sum_{i=1}^{\ell}T_{e_i}^{2}, \quad H=[E,F], \quad
\mathtt{E}(k)= \sum_{i=1}^{\ell}z_iT_{e_i}.
$$
Par un calcul analogue à celui de~\cite[Theorem~3.3]{H} on
montre que $(E,F,H)$ est un $\fsl(2)$-triplet d'\'el\'ements de
$\hec^W$ et que $H=\mathtt{E}(k) + g_k$ avec $g_k=\frac{\ell}{2} +
\sum_{\alpha \in R^+}k_\alpha r_\alpha$ (qui est un \'el\'ement
central de $\mathbb{C}W$).

\begin{Theoreme}
\label{thm41}
Soit $V$ un $\calH$-module irréductible de dimension finie. Il
existe un $\chi \in W\sphat\,$ tel que $V \simeq L(\chi)$.
\end{Theoreme}

\begin{dem}
  Décomposons le $W$-module $V$ en somme de composantes
  isotypiques : $V=\sum_{\chi \in W\sphat}V[\chi]$. Remarquons que
  si $P \in \calH^W$ on a $P.V[\chi] \subset V[\chi]$; en
  particulier, chaque $V[\chi]$ est un sous-$\fsl(2)$-module de
  dimension finie.  Soit $u$ de poids minimal $m \in - \N$ dans le
  $\fsl(2)$-module $V$.  \'Ecrivons $u = \sum_\chi u_\chi$ avec
  $u_\chi \in V[\chi]$. Comme $H$ laisse stable $V[\chi]$ pour
  tout $\chi$, chaque $u_\chi$ non nul est aussi un vecteur de
  poids minimal $m$ dans $V$.  On peut donc supposer que $u \in
  V[\chi]$ pour un $\chi \in W\sphat$.  De plus, pour tout $i \in
  \{1,..,\ell\}$ on a $[H,T_{e_i}]= -T_{e_i}$, donc
  $HT_{e_i}.u=(m-1)T_{e_i}.u$. Par minimalit\'e de $m$ il vient
  $T_{e_i}.u=0$, donc $\calS_+.u =0$.
  
  Posons $\CW.u = \boplus_j V_j$, $V_j \simeq V_\chi$. Soit $a
  =\sum_w a_w w \in \CW$ tel que $v=a.u \in V_j \smallsetminus
  \{0\}$. Puisque $T_{y}w=wT_{w^{-1}(y)}$ pour tout $y \in \fa$,
  il vient
  $$
  \textstyle{T_y.v = \sum_w a_w wT_{w^{-1}(y)}.u =0.}
  $$
  Par conséquent $V= \calH.v$ avec $\CW.v \simeq V_\chi$ et
  $\calS_+.v = 0$; le (f) de la Proposition~\ref{prop32} donne $V
  \simeq L(\chi)$.
\end{dem}

Rappelons que le module $M(\chi) = \calP\otimes \V$ est gradu\'e
par les $M_n(\chi)=\calP_n\otimes \V={\hec}_n.\v$, $n\geq0$.  En
posant $R_n(\chi) = R(\chi) \cap M_n(\chi)$ on en déduit la
graduation $R(\chi)= \boplus_{n\geq0}R_n(\chi)$ du radical de la
forme $(\phantom{.},\phantom{.})$ introduite au \S~\ref{sec3}. Le
quotient $L(\chi) = M(\chi)/R(\chi)$ est ainsi naturellement
gradué:
$$
L(\chi)= \boplus_{n\geq0}L_n(\chi) \ \, \text{avec
  $L_n(\chi)=M_n(\chi)/R_n(\chi)$.}
$$
Donc $L(\chi)$ est de dimension finie si, et seulement si, il
existe $n_0 \in \N$ tel que $L_n(\chi)=0$ pour $n \ge n_0$. Cette
condition équivaut à $M_n(\chi) = R_n(\chi)$, ou encore à:
$(\phantom{.},\phantom{.})$ est identiquement nulle sur
$M_n(\chi)$.


\section{\'Etude des $L(\chi)$ de dimension finie}
\label{sec5}
 
Nous allons donner une condition n\'ecessaire et suffisante,
portant sur $k$ et $\chi$, cf.~Remarque~\ref{rem55}, pour que
$L(\chi)$ soit de dimension finie et \'etudier plus pr\'ecisément
la structure d'un tel $L(\chi)$.  On conserve les notations des
sections~\ref{sec3} et~\ref{sec4}.

\begin{prop}
\label{prop51}
Soit $\chi \in W\sphat$.  Si $L_n(\chi)=0$, alors
$L_{n+1}(\chi)=0$.
\end{prop}

\begin{dem}
  Il s'agit de montrer que si $R_n(\chi)=M_n(\chi)$, alors
  $R_{n+1}(\chi)=M_{n+1}(\chi)$.  Soient $a,b \in \calP_{n+1}$ et
  $f,g\in {\mathbb{C}W}$, on doit montrer que
  $(a\otimes(f.\v),b\otimes(g.\v))=0$.  Il suffit de le faire
  lorsque $b$ est un mon\^ome de la forme $x_1x_2...x_{n+1}$, $x_j
  \in \aalg^{*}$; on \'ecrit $a=\sum_{x \in \aalg^{*}}xa_x$ avec
  $a_x \in \calP_n$.
  
  Soient $y \in \aalg$ et $w \in W$. De $[T_y,w] = w(T_{w^{-1}(y)}
  -T_y) = wT_{w^{-1}(y) -y}$ on tire que $[T_y,w] \in W\calS_{+}$.
  En \'ecrivant $[T_y,b]= \sum_{j} x_1\cdots [T_y,x_j] \cdots
  x_{n+1}$ et en utilisant la relation~1 du \S~\ref{sec1}, on
  montre que $[T_y,b] \in \calP_{n}\CW$. Alors, $T_ybw=bwT_y +
  [T_y,bw]= bwT_y + [T_y,b]w+b[T_y,w]$ et $\calS_+.\v = 0$
  impliquent $T_y bw.\v = [T_y,b]w.\v \in M_n(\chi)=\calP_n\otimes
  \V$.
  
  Calculons maintenant $(a\otimes(f.\v),b\otimes(g.\v))$. En
  utilisant la $\sigma$-sym\'etrie de $(\phantom{.},\phantom{.})$
  et la d\'efinition de $\sigma$ il vient
\begin{equation*}
\begin{align*}
  (a\otimes(f.\v),b\otimes(g.\v))
  &=\textstyle{\sum_{x\in\aalg^{*}}}
  (a_x\otimes(f.\v),\sigma(x)b\otimes(g.\v)) \\
  &=\textstyle{\sum_{x\in\aalg^{*}}}
  (a_x\otimes(f.\v),T_{B(x)}b\otimes(g.\v)).
\end{align*}
\end{equation*}
Il résulte du paragraphe précédent que pour tout $x \in \fa$:
$$
(a_x\otimes(f.\v),T_{B(x)}b\otimes(g.\v)) \in
(M_n(\chi),M_n(\chi)) = 0 \ \; \text{(par hypoth\`ese)}.
$$
Donc $(a\otimes(f.\v),b\otimes(g.\v))=0$.
\end{dem}

On en déduit:

\begin{cor}
\label{cor52}
Le module $L(\chi)$ est de dimension finie si, et seulement si, il
existe $n\in\mathbb{N}$ tel que $L_n(\chi)=0$.\qed
\end{cor}

Pour obtenir un critère plus pr\'ecis assurant la finitude de
$\dim L(\chi)$ nous allons utiliser le $\fsl(2)$-triplet $(E,F,H)$
d\'efini \`a la section pr\'ec\'edente.  On fera appel au
r\'esultat suivant \cite[Proposition~3.4]{H}:

\begin{Lemme}
\label{lem53}
Pour tout $p \in \calP_n$, on a
$$
\sigma(p)=(-1)^{n}\mathrm{ad}(F)^{n}(p)=
\textstyle{\sum_{j=0}^{n}}(-1)^{j}c_jF^{j}pF^{n-j}
$$
o\`u les $c_j$ sont des entiers ne d\'ependant que de $n$.
\qed
\end{Lemme}

\begin{Theoreme}
\label{thm54}
Soit $0 \ne v_\chi \in V_\chi$. Le module $L(\chi)$ est de
dimension finie si, et seulement si, il existe $m \geq 0$ tel que
$E^{m}.v_\chi \in R_{2m}(\chi)$.
\end{Theoreme}

\begin{dem} 
  Il est clair que la condition est n\'ecessaire, montrons qu'elle
  est suffisante.  Soit $n \in \mathbb{N}$. Tout \'el\'ement de
  $\hec_n.\v$ s'\'ecrit $ps.\v$ avec $p \in \calP_n$ et $s\in
  \mathbb{C}W$.  Soit $r\in\hec_n$ et calculons
  $(ps.\v,r.\v)=(p.\v,\sigma(s)r.\v)$.  Posons $t=\sigma(s)r
  \in\hec_n$.  On a $(p.\v,t.\v)=\pscal{\v}{\beta(p,t).\v}$ et
  d'apr\`es le lemme pr\'ec\'edent
  $\beta(p,t)=\pi(\sigma(p)t)=\sum_{j=0}^{n}(-1)^{j}
  c_j\pi(F^{j}pF^{n-j}t)$.  Or $F^{n-j}t \in
  \hec_{-2(n-j)+n}=\hec_{-(n-2j)} \subset \hec \calS^{+}$ pour
  $n-2j > 0$.  Donc $\pi(F^{j}pF^{n-j}t) \in \pi(\calH \calS_+) =
  0$ pour $n > 2j$. Compte tenu de $\sigma(E)=-F $ on obtient
  $$
  \textstyle{\beta(p,t)=\sum_{j=[\frac{n}{2}]}^{n}
    c_j\pi(\sigma(E^{j})pF^{n-j}t)=
    \sum_{j=[\frac{n}{2}]}^{n}c_j\beta(E^{j},pF^{n-j}t)}.
  $$
  (O\`u $[\phantom{.}]$ d\'esigne la partie enti\`ere.)  Par
  cons\'equent:
\begin{equation}\label{eq51}
\textstyle{(ps.\v,r.\v)=\sum_{j=[\frac{n}{2}]}^{n}c_j
  (E^{j}.\v,pF^{n-j}\sigma(s)r.\v)}.
\end{equation}

Supposons $E^{m}.v_\chi \in R_{2m}(\chi)$. Pour $j \ge m$ on a
donc
$$
(E^{j}.v_\chi,
M_{2j+2}(\chi))=(E^{m}.v_\chi,\sigma(E^{j-m})M_{2j+2}(\chi)) =0.
$$
L'\'equation~\eqref{eq51} fournit alors $(ps.\v,r.\v)= 0$ pour
tout $n\geq 2m$, c'est \`a dire $M_n(\chi) = R_n(\chi)$ pour $n
\ge 2m$, ce qui montre que $L(\chi)$ est de dimension finie.
\end{dem}

\begin{Remarque} \label{rem55}
  Puisque $E \in \calH^W$, on a $E^m.\v \in M(\chi)[\chi]$. La
  Remarque~\ref{rem34}(3) permet de pr\'eciser encore la condition
  obtenue dans le Th\'eor\`eme~\ref{thm54}: le module $L(\chi)$
  est de dimension finie si, et seulement si, il existe $m\geq 0$
  tel que $(E^{m}.\v,M_{2m}(\chi)[\chi])=0$.  Observons de plus
  que $\dim M_{2m}(\chi)[\chi] < \infty$ et que si $P \in
  M_{2m}(\chi)[\chi]$, $(E^m.v_\chi,P)$ est un polyn\^ome en les
  $k_\alpha$, $\alpha \in R^+$ (c'est en fait un polyn\^ome en les
  $k_i$ définis ci-dessous). Donc, si $\{P_1,\dots,P_s\}$ est une
  base de $M_{2m}(\chi)[\chi]$, la condition $E^m.v_\chi \in
  R_{2m}(\chi)$ équivaut à l'annulation des $s$ polyn\^omes
  $(E^m.v_\chi,P_j)$, $1 \le j \le s$.
\end{Remarque}

\subsection*{Le cas $\chi = \triv$} 
On a rappel\'e, cf.~Remarque~\ref{rem36}, que $M(\triv,k)=\calP$
et $R(\triv)=R(k)$. Gr\^ace \`a la Remarque~\ref{rem55}, la
condition du Th\'eor\`eme~\ref{thm54} se traduit par
$(E^{m},\calP^{W}_{2m})_k=0$, soit $F^{m}(\calP^{W}_{2m})=0$ pour
un $m \in \mathbb{N}$. Rappelons que $\calP^W
=\C[Q_1,\dots,Q_\ell]$ où les $Q_j$ sont des polyn\^omes homogènes
de degrés $d_1 \le d_2 \le \cdots \le d_\ell$ (appelés degrés
primitifs de $W$). La nullité de $F^{m}(\calP^{W}_{2m})$ est alors
équivalente à $F^m(Q_1^{a_1} \cdots Q_\ell^{a_\ell})=0$ pour tout
$(a_1,\dots,a_\ell) \in \N^\ell$ tel que $\sum_j a_j d_j = 2m$.

 \subsection*{Structure des $L(\chi)$ de dimension finie} Nous
 allons maintenant donner quel\-ques r\'esultats sur la forme du
 module $L(\chi)$ quand il est de dimension finie.  La structure
 de $\fsl(2)$-module de $L(\chi)$ lui donne une certaine
 sym\'etrie, du m\^eme type que celle observ\'ee en rang $1$,
 {cf.}~\cite[Theorem~9.2]{CM}.
 Soient $R_1,\dots,R_s$ les $W$-orbites dans $R$; on peut alors
 \'ecrire $R^+$ comme r\'eunion disjointe $\bigsqcup_{i=1}^s
 R_i^+$. On pose $k_\alpha = k_i$ pour $\alpha \in R_i^+$.
 Lorsque $R$ est irr\'eductible on a $s=2$ et l'on prend pour
 $R_1^{+}$, resp.~$R_2^{+}$, l'ensemble des racines courtes,
 resp.~longues (\'eventuellement vide), dans $R^{+}$.
 
 L'\'el\'ement $\sum_{\alpha\in R^{+}}k_\alpha r_\alpha$ \'etant
 central dans $\CW$, il op\`ere par multiplication par un scalaire
 sur tout $W$-module irr\'eductible $V$. Si $V$ est de type
 $\tau$, ce scalaire est
 $$
 a_\tau(k)= \sum_{i=1}^s k_i
 |R_i^{+}|\frac{\tau(r_i)}{\tau(1)},
 $$
 o\`u $\tau(r_i)$ d\'esigne la valeur commune des $\tau(r_\alpha)$
 pour $\alpha$ dans $R_i^{+}$.

\begin{prop}
\label{prop56}
{\rm 1.} Posons $b_\chi(k)= \frac{\ell}{2} + a_\chi(k)$.
L'\'el\'ement $H$ op\`ere sur $M_p(\chi)$ par le scalaire
$p+b_\chi(k)$.

\noindent {\rm 2.} Si $L(\chi)$ est de dimension finie
il existe un $m \in \N$ tel que:
\begin{enumerate}
\item[{\rm (i)}] $a_\chi(k) = -(m + \frac{\ell}{2}) < 0$;
\item[{\rm (ii)}] $L(\chi) = \boplus_{i=0}^{2m} L_i(\chi)$ et
  l'application $x \mapsto E^m.x$ est un isomorphisme de
  $W$-modules de $L_0(\chi) \simeq V_\chi$ sur $L_{2m}(\chi)$; si
  $x \in V_\chi \smallsetminus \{0\}$ on a
  $$
  m = \min\{p\in \N : E^{p+1}.x \in R_{2p+2}(\chi)\}.
  $$
\end{enumerate}
\end{prop}

\begin{dem}
  1. On rappelle, cf.~\S~\ref{sec4}, que $H=\mathtt{E}(k) + g_k$
  avec $g_k= \sum_{\alpha\in R^{+}} k_\alpha r_\alpha +
  \frac{\ell}{2}$.  Alors, par~\cite[Proposition~2.26]{DO},
  $\mathtt{E}(k)$ op\`ere sur $M_p(\chi)[\tau]$ par
  multiplication par le scalaire $p + a_\chi(k) -a_\tau(k)$.  Donc
  $H$ op\`ere sur $\boplus_{\tau\in W\sphat} \,
  M_p(\chi)[\tau]=M_p(\chi)$ par multiplication par
  $p+\frac{\ell}{2}+ a_\chi(k) = p +b_\chi(k)$.
  
  2. En passant au quotient modulo $R_p(\chi)$, on obtient que
  pour tous $p \geq 0$ et $x\in L_p(\chi)$, $H.x=(p+b_\chi(k))x$.
  Les \'el\'ements de poids minimal pour $H$ sont donc dans
  $L_0(\chi)= V_\chi$.  De plus, comme $[H,E]=2E$ et $[H,F]=-2F$,
  on a $H.(E.x)=(p+2+b_\chi(k))E.x$ et
  $H.(F.x)=(p-2+b_\chi(k))F.x$.
  
  Supposons maintenant $L(\chi)$ de dimension finie. On a donc
  $L(\chi)=\boplus_{i=0}^p L_i(\chi)$, avec $L_j(\chi)\neq 0$, $1
  \le j \le p$, cf.~Proposition~\ref{prop51}. Alors si $0 \ne x
  \in L_0(\chi)$, $x$ est vecteur propre de plus bas poids
  $b_\chi(k)$ pour l'action de $\fsl(2)$ et $U(\fsl(2)).x$ est un
  $\fsl(2)$-module irr\'eductible de dimension finie.  Il existe
  donc $m\in \mathbb{N}$ tel que $b_\chi(k)=-m$. D'o\`u (i).

\noindent (ii) De $U(\fsl(2)).x = \boplus_{i=0}^m \C E^i.x$ avec
$0 \ne E^m.x \in L_{2m}(\chi)$ vecteur de plus haut poids (\'egal
\`a $m$) on tire $2m \le p$. Soit $z \in L_p(\chi)$. Comme $E.z =
0$ et $H.z = (p-m)z$, le $\fsl(2)$-module $U(\fsl(2)).z$ est
simple de plus haut poids $p-m$. Il en d\'ecoule que $0 \ne
F^{p-m}.z \in L_{-2(p-m)+m}(\chi)$, donc $-2(p-m)+p=2m-p \geq 0$
et $p=2m$.

Comme $E^m \in \calH^W$, l'application (non nulle) $x \mapsto
E^m.x$ de $L_0(\chi) \simeq V_\chi$ vers $L_{2m}(\chi)$ est
injective. Il reste \`a montrer que son image est $L_{2m}(\chi)$.
Soit $y\in L_{2m}(\chi)$; c'est un vecteur de plus haut poids pour
$\fsl(2)$.  Par cons\'equent $E^{m}.(F^{m}.y)=y$, \`a une
constante pr\`es, avec $F^{m}.y \in L_0(\chi)$. La dernière
assertion découle immédiatement de la précédente.
\end{dem}

\begin{Remarques}
\label{rem57} 
\noindent (1) Supposons $L(\chi,k)$ de dimension finie et soit $m$
comme dans la Proposition~\ref{prop56}.  En raisonnant comme dans
la preuve du 2(ii) de la cette proposition, on peut montrer que
l'application $v \mapsto E^{n-1}v$ est un isomorphisme de
$W$-modules de $L_1(\chi)=M_1(\chi) / R_1(\chi)$ sur
$L_{2n-1}(\chi)$.
 
\noindent (2) Si 
$\chi(r_i)= 0$ pour tout $i$, il r\'esulte de la Proposition~\ref{prop56}(i)
que $L(\chi,k)$ est de dimension infinie. Cette condition est par exemple
v\'erifi\'ee lorsque $\chi=\chi\otimes\sgn$.

\end{Remarques}

\section{Exemples} 
\label{sec6}

{\em On suppose, dans toute cette section, que $R$
  irr\'eductible.} On a alors $\calP^W = \C[Q_1,\dots,Q_\ell]$
avec $d_1 = 2 < d_3 \le \dots \le d_\ell$ et l'on peut prendre
$Q_1=E$.

On rappelle que la multiplicit\'e $k$ est dite {\em singuli\`ere}
si $R(k)\neq 0$.  D'apr\`es \cite{DOJ}, pour $k$ constante, cela
\'equivaut \`a $k=\frac{j}{d_i} - p$ avec $1 \le i \le \ell$, $1
\le j \le d_i-1$ et $p \in \N^*$.  On dira que $k$ est {\em tr\`es
  singuli\`ere} si $L(\triv,k)$ est de dimension finie. (Rappelons
que c'est \'equivalent \`a l'existence d'un entier $p$ tel que
$F^{p}(\calP_{2p}^{W})=0$.)  Une multiplicit\'e tr\`es
singuli\`ere est évidemment singuli\`ere et l'on a vu que, si
$\tau$ est une repr\'esentation de dimension~$1$, cela équivaut à
$\dim L(\tau,k^\tau) < \infty$ (cf.~Remarque~\ref{rem36}).  Le but
de cette section est de donner des exemples de multiplicit\'es
tr\`es singuli\`eres pour certains syst\`emes de racines.

On pose $\hbar=b_\triv(k)=H(1)$, donc
$$
\hbar =-F(E)=k_1|R_1^{+}|+k_2|R_2^{+}|+ \frac{\ell}{2} = \half
|R^{+}|(k_1 + k_2) +\frac{\ell}{2}.
$$
Observons que $\hbar =k|R^+| + \frac{\ell}{2}$ lorsque toutes
les racines ont la m\^eme longueur.

\begin{Lemme}
\label{lem61}
Pour tout $p\in\N$, $F^{p+1}(E^{p+1})=(-1)^{p+1} (p+1)!
\prod_{i=0}^{p} (\hbar+i).$
\end{Lemme}

\begin{dem}
  En utilisant le fait que $H(P) = (\hbar + d)P$ pour $P \in
  \calP_d^W$ et la relation
\begin{equation} \label{eq61}
[F,E^s] = -HE^{s-1} + E[F,E^{s-1}]
\end{equation}
on montre par récurrence que $F(E^s)=-s(\hbar +s-1)E^{s-1}$, d'où
le lemme.
\end{dem}

On posera, si ces entiers existent,
$$
n = \min\{p\in \N : F^{p+1}(E^{p+1})=0\}, \quad m = \min\{s \in
\N: F^{s+1}(\calP_{2s+2}^W) = 0\}.
$$
Donc l'existence de $m$ équivaut à $\dim L(\triv,k) < \infty$
et l'on a alors $L(\triv,k) = \boplus_{i=0}^{2m} L_i(\triv,k)$,
cf.~\S~\ref{sec5}. (On pourra remarquer que si $L_1(\triv, k) \ne
0$, il est \'egal au $W$-module $\fa^*$.) L'existence de $n$
signifie que $\hbar = -n$, i.e.~$k = -\frac{2n+\ell}{|R|}$ lorsque
$k$ est constante.

\subsection*{Cas où $m=n$}
Lorsque $m=n$, cela d\'etermine les multiplicit\'es tr\`es
singuli\`eres. Ceci se produit lorsque $\hbar = -n$ et
$\calP_{2n+2} = \C E^{n+1}$. Donnons des exemples:
\begin{itemize}
\item Comme $\calP_2^{W}=\C E$ on a:
  $$
  \hbar=0 \, \iff \, L(\triv,k)=\C \iff R(k) = \calP_+.
  $$
  Si $k=k_1=k_2$, on obtient alors $k= -\frac{\ell}{|R|}$.
\item Si aucun degr\'e primitif n'est \'egal \`a $4$ (c'est par
  exemple le cas pour les types $\mathsf{A}_2$, $\mathsf{F}_4$,
  $\mathsf{E}_6$, $\mathsf{E}_7$, $\mathsf{E}_8$ et
  $\mathsf{G}_2$) on a $\calP_4^W=\C E^2$. Donc
  $$
  \hbar = -1 \, \iff \,
  L(\triv,k)=\boplus_{i=0}^{2}L_i(\triv,k).
  $$
\item De m\^eme, quand $d_j \ne 6$ pour tout $j$ (e.g.~pour le
  type $\mathsf{E}_8$) on obtient:
  $$
  \hbar=-2 \, \iff \, L(\triv,k)=\oplus_{i=0}^{4}L_i(\triv,k).
  $$
\item Supposons $R$ de type $\mathsf{A}_1$ et notons $\calP =
  \C[z]$. On a $\calP_{2n+2} = \C E^{n+1} = \C z^{2n+2}$, donc il
  découle des résultats de la section~\ref{sec5} que:
  $$
  \hbar = -n \, \iff \, L(\triv,k) = \C[z]/(z^{2n+1}).
  $$
  On retrouve ainsi ce qui a été obtenu dans~\cite{CM}: les
  multiplicit\'es singuli\`eres et tr\`es singuli\`eres
  co\"{\i}ncident et correspondent aux multiplicit\'es
  $k=-\frac{1}{2}-n$, $n \in \N$.
\end{itemize}

\subsection*{Type $\mathsf{A}_2$}
On a $\calP^{W}=\C[E,Q_2]$ o\`u $Q_2$ est un polyn\^ome homog\`ene
$W$-invariant de degr\'e trois.  Rappelons que l'ensemble des
valeurs singuli\`eres de $k$ est form\'e des $k$ {\em non entiers}
de la forme $-p/3$ ou $-p/2$, $p \in \N$.

On a $F(Q_2) \in \calP_{1}^W = 0$ et un calcul direct montre que
$F(Q_2^2) = \lambda E^2$ pour un $\lambda \in \C^*$ indépendant de
$k$. En rempla\c{c}ant $Q_2$ par $\sqrt{\lambda}Q_2$ on peut donc
écrire $\calP^W = \C[E,Q]$ avec $F(Q) = 0$, $F(Q^2) = E^2$.  On
pose
\begin{equation*} 
P_{n,r} = P_{n,r}(k) = F^n(E^{n-3r}Q^{2r}).
\end{equation*}                                    
Les $P_{n,r}$ sont des polyn\^omes en $k$, ou $\hbar$, et l'on a
$F^n(\calP_{2n}^W) = \sum_{i=0}^{[\frac{n}{3}]} \C P_{n,i}$.
Ainsi la multiplicit\'e $k$ est tr\`es singuli\`ere si, et
seulement si, il existe un entier positif ou nul $m$ tel que
$P_{m+1,r} = 0$ pour tout $r=0,\dots,[\frac{m+1}{3}]$.

On peut montrer le r\'esultat qui suit en utilisant
$[F,Q](Q^p)=F(Q^{p+1}) -QF(Q^p)$ et la formule~\eqref{eq61}.

\begin{Assertion} \label{ass62}
  Les polyn\^omes (en $\hbar$) $\Pnr$ v\'erifient la formule de
  r\'ecurrence suivante:
  $$
  \forall \, n,r \in \N, \quad \Pnpr = r(2r -1) \pnr{n}{r-1} -
  (n+1 -3r) (\hbar + n +3 r) \Pnr.
  $$
  En outre, $\Pnr = 0$ si $n < 3r$. \qed
\end{Assertion}   

Pour tous $n \in \N$ et $0 \le r \le [\frac{n}{3}]$ on d\'efinit
alors un polyn\^ome (en $\hbar$) de degr\'e $n-r$ par:
\begin{equation*}
\textstyle{f_{n,r}(\hbar) = \prod_{j=0}^{n-1}(\hbar +j) \Big/
\prod_{i=0}^{r-1}(\hbar +3i +2).
}
\end{equation*} 
Remarquons que $f_{n,r}(\hbar)$ s'annule pour $\hbar = -j$ avec $j
\in \{0,1,3, \dots,n-1\}$ et $j \not\equiv 2 \pmod{3}$. Donc, si
$n \not\equiv 2 \pmod{3}$ le polyn\^ome
$f_{n+1,\textstyle{[\frac{n+1}{3}]}}(\hbar)=
f_{n+1,\textstyle{[\frac{n}{3}]}}(\hbar)$ poss\`ede $-n$ pour
racine.  Observons \'egalement que $f_{n+1,r+1}(\hbar)$ divise
$f_{n+1,r}(\hbar)$.

Gr\^ace \`a la formule de r\'ecurrence obtenue dans
l'Assertion~\ref{ass62} on peut alors d\'emontrer:

\begin{Assertion} \label{ass63}
\noindent {\rm 1.} Soit $n,r \in \N$. Alors, si $a_{n,r} =  (-1)^{n+r}
\frac{n!}{3^r} \prod_{i=1}^r(2i-1)$, on a
$$
P_{n,r} = a_{n,r}f_{n,r}(\hbar).
$$

\noindent {\rm 2.} La multiplicit\'e $k$ est tr\`es singuli\`ere si et
seulement si $\hbar =-m \in - \N$ avec $m \not\equiv 2 \pmod{3}$.
\qed
\end{Assertion}
              
Compte tenu de l'assertion précédente, des
remarques~\ref{rem57}(2) et~\ref{rem36} nous pouvons énoncer:

\begin{Theoreme}
\label{thm64}
Soit $R$ un système de racines de type $\mathsf{A}_2$, $k
=k_\alpha$, $\alpha \in R$. Alors $\calH(k)$ possède des
représentations de dimension finie si, et seulement si, il existe
un $m \in \N$ non congru à $2$ modulo $3$ tel que $3k + 1 = \pm
m$. De plus, $3k + 1 = -m$ $\iff$ $\dim L(\triv,k) < \infty$, et
$3k + 1 = m$ $\iff$ $\dim L(\sgn,k) < \infty$. \qed
\end{Theoreme}

\subsection*{Type $\mathsf{B}_2$} Rappelons que  $k_1$,
resp.~$k_2$, d\'esigne la valeur de $k$ sur les racines courtes,
resp.~longues; on a donc $\hbar = 2(k_1 + k_2) + 1$. \'Ecrivons
$\calP= \C[x_1,x_2]$ de sorte que $R^+=\{x_1,x_2,x_1 \pm x_2\}$ et
$W$ soit engendr\'e par $r=r_{x_1-x_2}$, $s_i=r_{x_i}$, $i=1,2$.
On a alors $\calP^W=\C[E,Q]$ o\`u $Q=x_1^2x_2^2 \in \calP_4^W$.

Soit $n \in \N$. Une base de $\calP_{2n+2}^W$ est $\{E^{n+1 -
  2r}Q^r$, $r = 0,\dots,[\frac{n+1}{2}]\}$.  On d\'efinit des
polyn\^omes en $k_1,\hbar$ (i.e.~$k_1,k_2$) par
\begin{equation*} 
P_{n,r} = P_{n,r}(k_1,\hbar) = F^n(E^{n-2r}Q^{r}).
\end{equation*}
La condition $F^{n+1}(\calP_{2n+2}^W)=0$ est donc \'equivalente
\`a $P_{n+1,r} = 0$ pour $r=0,\dots,[\frac{n+1}{2}]$.  On obtient
ais\'ement le r\'esultat suivant:

\begin{Assertion}
\label{ass64}
Les polyn\^omes $\Pnr$ v\'erifient la formule de r\'ecurrence:
$$
\forall \, n,r \in \N, \quad \Pnpr = -2r (2k_1+2r -1)
\pnr{n}{r-1} - (n+1 -2r) (\hbar + n + 2r) \Pnr.
$$
(Avec la convention $P_{n,-1} = 0$.) \qed
\end{Assertion}

Pour tous $n \in \N$ et $0 \le r \le [\frac{n}{2}]$, on pose:
\begin{equation*}
\textstyle{g_{n,r}(k_1,\hbar) = \prod_{i=1}^r (2k_1 + 2i
  -1)\prod_{j=0}^{n-1}(\hbar +j) \Big/  
\prod_{i=1}^{r}(\hbar +2i -1). 
}
\end{equation*}

Une r\'ecurrence et l'assertion~\ref{ass64} donnent:

\begin{Assertion}
\label{ass65}
Soient $n,r \in \N$. Alors, $\Pnr = (-1)^n n! \,
g_{n,r}(k_1,\hbar)$.  \qed
\end{Assertion}

Rappelons que $W\sphat =\{\triv, \sgn, \mathrm{std},
\chi_1,\chi_2\}$, o\`u $\mathrm{std}$ est la repr\'esentation
standard $W \hookrightarrow \mathrm{GL}(\fa)$, $\chi_1$, $\chi_2$
sont de dimension $1$ donn\'ees par $\chi_1(r) = 1$, $\chi_1(s_1)
= \chi_1(s_2)=-1$, $\chi_2 = \chi_1 \otimes \sgn$.  En utilisant
la formule pr\'ec\'edente, les remarques~\ref{rem57}(2)
et~\ref{rem36}, on obtient le th\'eor\`eme qui suit. Il fournit
les conditions n\'ecessaires et suffisantes pour que $\calH(k)$
poss\`ede des repr\'esentations de dimension finie.

\begin{Theoreme}
\label{thm66}
Soit $R$ un système de racines de type $\mathsf{B}_2$.

\noindent {\rm 1.} La condition  $\dim L(\triv,k) < \infty$
\'equivaut \`a l'une des deux conditions suivantes:
\begin{enumerate}
\item[{\rm (a)}] $\hbar = -n$ pour un $n \in \N$ pair;
\item[{\rm (b)}] $\hbar = -n$ et $k_1 = -\frac{2j-1}{2}$ avec $j
  \in \{1,\dots,\frac{n+1}{2}\}$ pour un $n \in \N$ impair.
\end{enumerate}

\noindent {\rm 2.} Le module  $L(\mathrm{std},k)$ est de
dimension infinie.

\noindent {\rm 3.} Si $\tau \in \{\sgn, \chi_1, \chi_2\}$ on a:
$\dim L(\tau,k) < \infty$ $\iff$ $\dim L(\triv,k^\tau) < \infty$.
\qed
\end{Theoreme}

\subsection*{Type $\mathsf{G}_2$} Soit $e_1, e_2$ une base
orthonorm\'ee de $\aalg_{\R}$, de base duale $x_1, x_2$.  On pose:
$z=x_1+ix_2$, $\bar{z}=x_1-ix_2$ et $T= \half (T_{e_1} - i
T_{e_2})$, $ \bar{T} = \half (T_{e_1}+ i T_{e_2})$.  Ainsi $2E = z
\bar{z}$ et $F = - 2 T \bar{T}$.  Soit $Q^{'}=z^6+\bar{z}^6$.
Comme $W$ est le groupe di\'edral d'ordre $12$, on a
$\hbar=1+3(k_1+k_2)$, $\calP^W=\C[E,Q^{'}]$ et (cf.~\cite{DOJ})
la multiplicit\'e $k=(k_1,k_2)$ est singuli\`ere si, et seulement
si, il existe $n \in \N$ tel que l'on soit dans l'un des cas
suivants:
\begin{enumerate}
\item[(i)] $k_1 =-\frac{1}{2} -n$, $k_2 \in \C$;
\item[(ii)] $k_2=-\frac{1}{2} -n$, $k_1 \in \C$;
\item[(iii)] $3(k_1+k_2)= -j-3n$ o\`u $j \in \{1,2,4,5\}$.
\end{enumerate}

Posons $\kappa=k_2-k_1$; il existe alors une constante $\lambda$
(non-nulle et ind\'ependante de $k_1,k_2$) telle que $F(\lambda Q^{'})=\kappa
E^2$. On pose $Q=\lambda Q^{'}$ et, pour tous $n\in \N$ et $r \in \{0, \dots,
[\frac{n}{3}]\}$:
$$
P_{n,r}=F^n(Q^rE^{n-3r}).
$$
Ces polyn\^omes en $k_1,k_2$ (ou $\kappa,\hbar$) engendrent
$F^n(\calP_{2n}^W)$. On montre par r\'ecurrence le r\'esultat
suivant.


\begin{Assertion} \label{ass68}
  Pour tous $(n,r) \in \N^2$, on a:
\begin{equation*}
P_{n+1,r} = -(n+1 -3r) (\hbar +n + 3r) P_{n,r} + r\kappa
P_{n,r-1} - \frac{r(r-1)}{9} P_{n,r-2}.
\end{equation*}
(Avec la convention que $P_{n,r}= 0$ pour $r < 0$ ou $n<3r$.)
\qed
\end{Assertion}
 
On obtient  une suite de polyn\^omes (en
$\kappa,\hbar$), $\{\Phi_p(\hbar,\kappa)\}_{p\in\N}$, en posant
$\Phi_0=1$, $\Phi_1=\kappa$ et pour tout $p>1$:
$$
\Phi_p(\hbar,\kappa)=\kappa\Phi_{p-1}(\hbar,\kappa)+\frac{p-1}{3}
(\hbar+3p-4) \Phi_{p-2}(\hbar,\kappa).
$$
D\'efinissons aussi, pour $n\in\N$ et $r\in\{0,\dots,
[\frac{n}{3}]\}$, les polyn\^omes:
$$
a_{n,r}(\hbar) = \prod_{i=0}^{n-1}(\hbar+i) \Big/
\prod_{j=0}^{r-1}(\hbar +2+3j).
$$
On v\'erifie alors que pour $n\in\N$ et $r\in\{0,\dots,
[\frac{n}{3}]\}$ on a:
$$
P_{n,r}(\hbar,\kappa)=(-1)^{n+r}\frac{n!}{3^r}\Phi_r(\hbar,\kappa)
a_{n,r}(h).
$$ 
Ainsi l'existence d'un entier $n\in\N$ minimal tel que
$F^n(\calP_{2n}^W)=0$ est \'equivalente \`a la r\'ealisation de
l'une des deux conditions suivantes :
\begin{enumerate} 
\item[(a)] si $n \not\equiv 0 \pmod{3}$, alors $\hbar=-(n-1)$;
\item[(b)] si $n=3r$, alors $\phi_r(\kappa)= \Phi_r(\kappa,
  -(3r-1))=0$. 
\end{enumerate}
Quelques calculs fournissent la conjecture suivante (v\'erifi\'ee
jusqu'\`a $2q+1=31$):
\begin{equation}
\label{eq62}
\phi_{2q}(\kappa)=\prod_{j=1}^{q}(\kappa^2-(2j-1)^2) , \quad
\phi_{2q+1}(\kappa)=\kappa \prod_{j=1}^{q}(\kappa^2-(2j)^2).
\end{equation}

Supposons que~\eqref{eq62} soit v\'erifi\'ee.
La multiplicit\'e $k$ est alors tr\`es singuli\`ere si, et
seulement si, il existe $n\in \N$ tel que:
\begin{enumerate} 
\item[(a)] soit $n \not\equiv 0 \pmod{3}$ et $\hbar=-(n-1)$;
\item[(b)] soit $n \equiv 0 \pmod{3}$ et si $\frac{n}{3}$ est pair
  (resp.~impair), alors $|\kappa|$ est un entier impair
  (resp.~pair) inf\'erieur \`a $\frac{n}{3}$. 
\end{enumerate}  

Soit $\alpha_1$, resp.~$\alpha_2$, une racine courte, 
resp.~longue. On a $W\sphat =\{\triv, \sgn, \tau, \sgn\otimes\tau, 
\mathrm{std}, \mathrm{std}\otimes\tau\}$, o\`u $\tau$ est 
la representation irr\'eductible de dimension $1$ donn\'ee par 
$\tau(r_{\alpha_1})=1$ et $\tau(r_{\alpha_2})=-1$. En utilisant
la caract\'erisation pr\'ec\'edente, les remarques~\ref{rem57}(2)
et~\ref{rem36}, on obtient le th\'eor\`eme suivant. Il 
d\'etermine les multiplicit\'es pour lesquelles $\calH(k)$
admet des repr\'esentations de dimension finie.

\begin{Theoreme}
\label{thm69}
Soit $R$ un syst\`eme de racines de type $\mathsf{G}_2$. 

\noindent {\rm 1.} On suppose \eqref{eq62} v\'erifi\'ee. Alors la
multiplicit\'e $k=(k_1,k_2)$ est tr\`es singuli\`ere, i.e.~$\dim
L(\triv,k) < \infty$, si et seulement si il existe $n \in \N^*$ tel
que
\begin{enumerate}
\item[{\rm (a)}]  $k_1 + k_2 =
  -\frac{n}{3}$ si $n \not\equiv 0 \pmod{3}$; 
\item[{\rm (b)}] $(k_1,k_2)$ ou
  $(k_2,k_1)$ est de la forme $(-\half(\frac{n}{3}+2j-1),
  -\half(\frac{n}{3}- 2j +1))$, $1 \le j \le \frac{n}{6}$, si $n
  \equiv 0 \pmod{6}$; 
\item[{\rm (c)}]  $(k_1,k_2)$ ou
  $(k_2,k_1)$ est de la forme $(-\half(\frac{n}{3}+2j),
  -\half(\frac{n}{3}- 2j))$, $0 \le j \le \frac{n-3}{6}$, si $n
  \equiv 3 \pmod{6}$. 
\end{enumerate}

\noindent {\rm 2.} Si $\chi \in \{\sgn,\tau, \sgn \otimes \tau\}$,
alors $\dim L(\chi,k) < \infty$ si et seulement si la
multiplicit\'e $k^\chi$ est tr\`es singuli\`ere.

\noindent {\rm 3.} Si $\chi \in \{\mathrm{std},
\tau\otimes\mathrm{std}\}$, alors  $\dim L(\chi,k) = \infty$. 
\qed
\end{Theoreme}


\end{document}